
\newcount\secno
\newcount\prmno
\def\section#1{\vskip1truecm
               \global\def\currenvir{section}
               \global\advance\secno by1\global\prmno=0
               {\bf \number\secno. {#1}}
               \smallskip}

\def\subsection{\global\def\currenvir{subsection}
                \global\advance\prmno by1
               \smallskip  \ind{ (\number\secno.\number\prmno) }}
\def\subsec{\global\def\currenvir{subsection}
                \global\advance\prmno by1\smallskip
                { (\number\secno.\number\prmno)\ }}

\def\proclaim#1{\global\advance\prmno by 1
                {\bf #1 \the\secno.\the\prmno$.-$ }}

\long\def\th#1 \enonce#2\endth{%
   \medbreak\proclaim{#1}{\it #2}\global\def\currenvir{th}\smallskip}

\def\rem#1{\global\advance\prmno by 1
{\it #1} \the\secno.\the\prmno$.-$ }

\magnification 1250 \pretolerance=500 \tolerance=1000
\brokenpenalty=5000 \mathcode`A="7041 \mathcode`B="7042
\mathcode`C="7043 \mathcode`D="7044 \mathcode`E="7045
\mathcode`F="7046 \mathcode`G="7047 \mathcode`H="7048
\mathcode`I="7049 \mathcode`J="704A \mathcode`K="704B
\mathcode`L="704C \mathcode`M="704D \mathcode`N="704E
\mathcode`O="704F \mathcode`P="7050 \mathcode`Q="7051
\mathcode`R="7052 \mathcode`S="7053 \mathcode`T="7054
\mathcode`U="7055 \mathcode`V="7056 \mathcode`W="7057
\mathcode`X="7058 \mathcode`Y="7059 \mathcode`Z="705A
\def\spacedmath#1{\def\packedmath##1${\bgroup\mathsurround =0pt##1\egroup$}
\mathsurround#1
\everymath={\packedmath}\everydisplay={\mathsurround=0pt}}
 \spacedmath{2pt}

\def\iso{\vbox{\hbox to .8cm{\hfill{$\scriptstyle\sim$}\hfill}
\nointerlineskip\hbox to .8cm{{\hfill$\longrightarrow $\hfill}} }}
\def\sdir_#1^#2{\mathrel{\mathop{\kern0pt\oplus}\limits_{#1}^{#2}}}

\font\eightrm=cmr8 \font\sixrm=cmr6

\def\pc#1{\tenrm#1\sevenrm}
\def\tx{\kern-1.5pt -}
\def\cqfd{\kern 2truemm\unskip\penalty 500\vrule height 4pt depth 0pt width
4pt\medbreak} 
\def\no{n\up{o}\kern 2pt}
\def\ind{\par\hskip 1truecm\relax}

\font\pal=cmsy7

\def\sp#1{{\cal S}\kern-1pt\raise-1pt\hbox{\pal P}^{}_C(#1)}

\frenchspacing
\input xy
\xyoption{all}
\input amssym.def
\input amssym
\vsize = 25truecm \hsize = 16.1truecm \voffset = -.5truecm
\parindent=0cm
\baselineskip15pt \overfullrule=0pt

\vglue 2.5truecm \font\Bbb=msbm10

\centerline{\bf Norm groups and class fields of formally real
quasilocal fields}
\bigskip

\centerline{I.D. Chipchakov\footnote{$^{\ast}$}{Partially
supported by Grant MI-1503/2005 of the Bulgarian National Science Fund.}}
\par
\vskip1truecm Abstract. This paper establishes a relationship
between finite extensions and norm groups of formally real
quasilocal fields, which yields a generally nonabelian local class 
field theory, including analogues to the fundamental correspondence,
the local reciprocity law and the norm limitation theorem.
\par
\vskip1.truecm
\centerline{{\bf 1. Introduction}}
\par
\medskip
Let $E$ be a field, $\overline E$ an algebraic closure of $E$,
Fe$(E)$ the set of finite extensions of $E$ in $\overline E$ and
Nr$(E) = \{N(R/E)\colon \ R \in {\rm Fe}(E)\}$ the set of norm
groups of $E$. The purpose of this paper is to characterize the
elements of Nr$(E)$ among the subgroups of the multiplicative group
$E ^{\ast }$ of $E$, under the hypothesis that $E$ is formally real
and quasilocal. We show that the index $i(L/E)$ of $N(L/E)$ in $E
^{\ast }$ is at most equal to the degree $[L\colon E]$, for each $R
\in {\rm Fe}(E)$; when $i(L/E) = [L\colon E]$, $L$ is called a class
field of $E$. Our main results prove, for each $R \in {\rm Fe}(E)$,
the existence of class fields $\widetilde R$ of $E$ with
$N(\widetilde R/E) = N(R/E)$; as it turns out, these $\widetilde R$
are isomorphic over $E$ to a uniquely determined intermediate field
cl$(R/E)$ of $R/E$. We also show that the set $I({\rm cl}(R/E)/E)$
of extensions of $E$ in cl$(R/E)$ consists of class fields and the
natural map of $I({\rm cl}(R/E)/E) \to {\rm Nr}(E)$ is injective
with an image equal to the set of subgroups of $E ^{\ast }$
including $N(R/E)$. This allows us to describe Nr$(E)$ and to
characterize class fields of $E$, up-to an $E$-isomorphism. The form
of the main results is particularly simple when $E$ is strictly
quasilocal (abbr, SQL), i.e. finite extensions of $E$ are strictly
PQL-fields. Then Fe$(E)$ consists of class fields of $E$, and the
union $\Sigma _{1} \cup \Sigma _{2}$ is a system of representatives
of the $E$-isomorphism classes of Fe$(E)$, where $\Sigma _{0} =
\{\Phi _{0} \in {\rm Fe}(E)\colon \ \Phi _{0} \subseteq E _{\rm
rc}$, for a fixed real closure $E _{\rm rc}$ of $E$ in $\overline
E$, and $\Sigma _{1} = \{\Phi _{1} \in {\rm Fe}(E)\colon \ \sqrt{-1}
\in \Phi _{1}\}$. Also, the canonical correspondence $\Sigma _{0}
\cup \Sigma _{1} \to {\rm Nr}(E)$ is bijective and maps compositums
into intersections and intersections into subgroup products.
\par
\medskip
As in [8], the basic field-theoretic notions needed for describing
the main results of this paper are the same as those used, e.g. in
[12, 14, 15] and [16]. As usual, $E _{\rm sep}$ denotes the
separable closure of a field $E$ in $\overline E$, and $E ^{\ast n}
= \{\alpha ^{n}\colon \ \alpha \in E ^{\ast }\}$, for each $n \in
\hbox{\Bbb N}$. The field $E$ is called formally real, if $-1$ is
not presentable over $E$ as a finite sum of squares; we say that $E$
is nonreal, otherwise. In the former case, $E$ is said to be
Pythagorean, if $E ^{\ast 2}$ is additively closed. Throughout this
paper, $\overline P$ denotes the set of prime numbers, Br$(E) _{p}$
is the $p$-component of the Brauer group Br$(E)$, for each $p \in
\overline P$, $\Pi (E)$ is the set of all $\pi \in \overline P$, for
which the absolute Galois group ${\cal G} _{E} = {\cal G}(E _{\rm sep}/E)$ is of
nonzero cohomological $\pi $-dimension cd$_{\pi } ({\cal G} _{E})$, $\Pi
_{m} (E) = \{\pi _{m} \in \overline P\colon \ {\rm cd}_{\pi _{1}} ({\cal G}
_{E}) = m\}$, for each $m \in \hbox{\Bbb N}$, and $P(E)$ stands for
the set of those $p \in \overline P$, for which $E$ is properly
included in its maximal $p$-extension $E (p)$ in $E _{\rm sep}$. In
what follows, the considered division algebras are supposed to be
associative with a unit, homomorphisms of profinite groups are
assumed to be continuous, and Galois groups are viewed as profinite
with respect to the Krull topology. For convenience of the reader,
we recall that $E$ is said to be PQL, if every cyclic extension $F$
of $E$ is embeddable as an $E$-subalgebra in each finite-dimensional
central division $E$-algebra $D$ of Schur index ind$(D)$ equal to
$[F\colon E]$. The field $E$ is called quasilocal, if its finite
extensions are PQL. We say that $E$ is a strictly PQL-field, if it
is PQL and Br$(E) _{p} \neq \{0\}$, for every $p \in P(E)$. Let us
note that PQL-fields and quasilocal fields appear naturally in the
process of characterizing basic types of stable fields with
Henselian valuations (see [4; 7] and the references there). Recall
also that strictly PQL-fields admit one-dimensional local class
field theory (abbr, LCFT), and SQL-fields are those whose finite
extensions have such a theory (see Section 2 and [8, Theorem 1.1 and
Corollary 3.6]).
\par
\medskip
The motivation and the aims of this research are determined by the
role of orderings in field theory and by their place in the study of
the PQL-property (see [9, Theorem 1.2 and Proposition 6.4], [7,
Theorem 3.1 and Lemma 3.5] and [4, Sect. 3]). As in the case of
SQL-fields with Henselian discrete valuations [5], our main results
yield a generally nonabelian one-dimensional LCFT (for a
presentation of the classical LCFT, see [11]). However, they
describe not only the groups in Nr$(E)$ but also the subgroups of $E
^{\ast }$ of finite indices. The first one of them can be stated as
follows:
\par
\medskip
{\bf Theorem 1.1.} {\it Let $E$ be a formally real quasilocal field.
Then:}
\par
(i) {\it Each $U \in {\rm Nr}(E)$ equals $N(\widetilde U/E)$, for
some class field $\widetilde U$ of $E$, which is uniquely determined
by $U$, up-to an $E$-isomorphism;}
\par
(ii) {\it A class field {\rm cl}$(U)$ of $E$ corresponding to a
group $U \in {\rm Nr}(E)$ embeds as an $E$-subalgebra in a field
$\Psi \in {\rm Fe}(E)$ if and only if $N(\Psi /E) \subseteq U$; if
$N(\Psi /E) = U$, then the $E$-isomorphic copy of {\rm cl}$(U)$ in
$\Psi $ is unique;}
\par
(iii) Nr$(E)$ {\it equals the set of subgroups of $E ^{\ast }$ of
finite indices not divisible by any $\pi \in \Pi _{1} (E) $, and
class fields of $E$ are precisely its finite extensions of
degrees with the same arithmetic property.}
\par
\medskip
Before stating the second main result, note that if $E$ is a
formally real quasilocal field, then ${\cal G} _{E}$ is metabelian (see
Section 3), so it possesses a unique conjugacy class of (closed)
Hall pro-$\Pi $-subgroups, for each subset $\Pi \subseteq \overline
P$.
\par
\medskip
{\bf Theorem 1.2.} {\it In the setting of Theorem 1.1, let $\Lambda
_{1}$ be the fixed field of some Hall pro-$\Pi _{1} (E)$-subgroup $H
_{\Pi _{1} (E)}$ of ${\cal G} _{E}$, and let $\Lambda _{0}$ be a real
closure of $E$ in $\Lambda _{1}$. Assume also that $\Sigma _{0} =
\{\Theta _{0} \in {\rm Fe}(E)\colon \ \Theta _{0} \subseteq \Lambda
_{0}\}$, and $\Sigma _{1} = \{\Theta _{1} \in {\rm Fe}(E)\colon \
\Theta _{1} \subseteq \Lambda _{1},$ $2 \vert [\Theta _{1}\colon
E]\}$. Then $\Lambda _{0}$ is uniquely determined, up-to an
$E$-isomorphism, $\Lambda _{1}/E$ is normal and the following
assertions are true:}
\par
(i) $\Sigma _{0} \cup \Sigma _{1}$ {\it is a system of
representatives of the class fields of $E$; the correspondence of
$\Sigma _{0} \cup \Sigma _{1}$ into {\rm Nr}$(E)$, defined by the
rule $\Phi \to N(\Phi /E)$, is bijective and satisfies the
conditions $N(L _{1}L _{2}/E) = N(L _{1}/E) \cap N(L
_{2}/E)$ and $N((L _{1} \cap L _{2})/E) = N(L _{1}/E).$ $N(L
_{2}/E)$, whenever $L _{j} \in \Sigma _{0} \cup \Sigma _{1}$, $j
= 1, 2$.}
\par
(ii) {\it For each $n \in \hbox{\Bbb N}$, $E ^{\ast n} = E ^{\ast
n(E)}$ and $E ^{\ast n}$ is of index $n(E)n(E) _{1}$ in $E ^{\ast
}$, where $n(E)$ and $n(E) _{1}$ are the greatest integers, such
that $n(E) _{1} \vert n(E) \vert n$, $4 \not\vert n(E)$, $2 \not\vert
n(E) _{1}$, $\pi \not\vert n(E)$ and $\pi _{1} \not\vert n(E) _{1}$,
for any $\pi \not\in \Pi (E)$ and $\pi _{1} \in \Pi _{1} (E)$. More
precisely, the quotient group $E ^{\ast }/E ^{\ast n}$ is isomorphic
to the direct product of the cyclic groups $C _{n(E)}$ and $C _{n(E)
_{1}}$.}
\par
(iii) $\Sigma _{1}$ {\it coincides with the set of finite Galois
extensions of $E$ in $\Lambda _{1}$, and $\sqrt{-1} \in M$, for
every $M \in \Sigma _{1}$.}
\par
\medskip
Let $E$ be a formally real quasilocal field and $D(E)$ the maximal
divisible subgroup of $E ^{\ast }$. Then $D(E) = \cap _{m=1}
^{\infty } E ^{\ast m}$, and the group $E ^{\ast }/D(E)$ can be 
endowed with a uniquely determined structure of a totally 
disconnected topological group so that the subgroups of $E ^{\ast
}/D(E)$ of finite indices form a system of neighbourhoods of unity
(cf. [16, Ch. 6, Theorem 9]). Also, Theorem 1.1 (iii) shows that the
groups $U/D(E)\colon \ U \in {\rm Nr}(E)$, are open in $E ^{\ast
}/D(E)$ with respect to that topology. It proves that the converse
is true if and only if $\Pi _{1} (E)$ is empty (or equivalently, $E$
is SQL, see Remark 3.2 (ii)). For a similar result in the nonreal
case, we refer the reader to [9, (6.1)].
\par
\medskip
Theorems 1.1 and 1.2 are proved in Section 3, and Section 2 includes preliminaries needed for the purpose. Our proof relies on the fact
(see [7, Lemma 3.5] and [3, (3.3)]) that a formally real field $E$
is hereditarily Pythagorean (i.e. its formally real extensions in
$\overline E$ are Pythagorean) with a unique ordering if and only if
all $F \in {\rm Fe}(E)$ are $2$-quasilocal fields in the sense of
[7].
\par
\vskip1.9truecm
\centerline {\bf 2. Preliminaries}
\par
\medskip
Let $E$ be a field and $\Omega (E)$ the set of finite abelian
extensions of $E$ in $E _{\rm sep}$. We say that $E$ admits
(one-dimensional) LCFT, if the mapping of $\Omega (E)$ into
Nr$(E)$ defined by the rule $F \to N(F/E)\colon \ F \in \Omega
(E)$, is injective and satisfies the following condition:
\par
(2.1) For each $M _{1}, M _{2} \in \Omega (E)$, $N(M _{1}M _{2}/E) =
N(M _{1}/E) \cap N(M _{2}/E)$ and
\par \noindent
$N((M _{1} \cap M _{2})/E) = N(M _{1}/E)N(M _{2}/E)$.
\par
Our approach to the study of fields with LCFT is based on the
following two lemmas (proved, e.g. in [6]).
\par
\medskip
{\bf Lemma 2.1.} {\it Let $E$ be a field and let $L$, $L _{1}$ and 
$L _{2}$ be elements of {\rm Fe}$(E)$, such that $L = L _{1}L _{2}$ 
and {\rm g.c.d.}$([L _{1}\colon E], [L _{2}\colon E]) = 1$. Then 
$N(L/E) = N(L _{1}/E) \cap N(L _{2}/E)$,
\par \noindent
$N(L _{1}/E) = E ^{\ast } \cap N(L/L _{2})$, and there is a group 
isomorphism 
$E ^{\ast }/N(L/E) \cong $
\par \noindent
$E ^{\ast }/N(L _{1}/E) \times E ^{\ast }/N(L _{2}/E)$.}
\par
\medskip
{\bf Lemma 2.2.} {\it Let $E$ and $M$ be fields such that $M \in
\Omega (E)$ and $M \neq E$. Let also $P(M/E) = \{p \in \overline
P\colon \ p \vert [M\colon E]\}$, and $M _{p} = M \cap E (p)$, for
each $p \in P(M/E)$. Then $N(M/E) = \cap _{p \in P(M/E)} N(M
_{p}/E)$ and $E ^{\ast }/N(M/E)$ is isomorphic to the direct
product of the groups $E ^{\ast }/N(M _{p}/E)\colon \ p \in P(M/E)$.}
\par
\medskip
The following lemma is known (it is a special case of [17, Ch. II,
Proposition 4]) but we include its direct proof due to its brevity
and simplicity.
\par
\medskip
{\bf Lemma 2.3.} {\it For a field $E$ and a given $p \in \overline
P$, the following conditions are equivalent:}
\par
(i) Br$(E ^{\prime }) _{p} = \{0\}${\it , for every extension $E
^{\prime }$ of $E$ in $\overline E$;}
\par
(ii) {\it The exponent $e _{1}$ of the group $E _{1} ^{\ast }/N(E
_{2}/E _{1})$ is not divisible by $p$, for any pair $(E _{1}, E
_{2})$ of finite extensions of $E$ in $E _{\rm sep}$, such that $E
_{1} \subseteq E _{2}$.}
\par
\medskip
{\it Proof.} (i)$\to $(ii): Let $E _{2} ^{\prime } \subseteq E _{\rm
sep}$ be the normal closure of $E _{2}$ over $E _{1}$, $e _{1}
^{\prime }$ the exponent of $E _{1} ^{\ast }/N(E _{2} ^{\prime }/E
_{1})$, and $E _{1} ^{\prime }$ the fixed field of a Sylow
$p$-subgroup of ${\cal G}(E _{2} ^{\prime }/E _{1})$. By [7, Lemma 4.2
(ii)], $N(E _{2} ^{\prime }/E _{1} ^{\prime }) = E _{1} ^{\prime
\ast }$, whence $N(E _{2} ^{\prime }/E _{1}) = N(E _{1} ^{\prime }/E
_{1})$ and $e _{1} ^{\prime } \vert [E _{1} ^{\prime }\colon E
_{1}]$. Since $N(E _{2} ^{\prime }/E _{1}) \subseteq N(E _{2}/E
_{1})$, $E _{1} ^{\ast }/N(E _{2}/E _{1})$ is a quotient of $E _{1}
^{\ast }$ $/N(E _{2} ^{\prime }/E _{1})$, so we have $e _{1} \vert e
_{1} ^{\prime } \vert [E _{1} ^{\prime }\colon E _{1}]$, proving
that (i)$\to $(ii).
\par
(ii)$\to $(i): Suppose that Br$(F) _{p} \neq \{0\}$, for some
$F$ of $E$ in $\overline E$. Then there exists $F _{0} \in {\rm
Fe}(E)$ with $F _{0} \subseteq F \cap E _{\rm sep}$ and Br$(F
_{0}) _{p} \neq \{0\}$ (cf. [17, Ch. II, 2.2 and 2.3] and [7,
(1.3)]). Hence, by [13, Sect. 4, Theorem 2], there is a central
division $F _{0}$-algebra $\Delta _{0}$ of index $p$. Moreover, by
[15, Sect. 15.2], $E _{\rm sep}$ has a subfield $F _{1} \in {\rm
Fe}(F _{0})$, such that $p \not\vert [F _{1}\colon F _{0}]$ and
$\Delta _{0} \otimes _{F _{0}} F _{1}$ is a cyclic division $F
_{1}$-algebra. This means that there is a cyclic extension $F _{2}$
of $F _{1}$ of degree $p$, embeddable in $\Delta _{0} \otimes _{F
_{0}} F _{1}$ over $F _{1}$. Therefore, by [15, Sect. 15.1,
Proposition b], $N(F _{2}/F _{1}) \neq F _{1} ^{\ast }$, which
completes the proof of Lemma 2.3.
\par
\medskip
For convenience of the reader, we give a proof of the following
lemma, which presents results on hereditarily Pythagorean fields
(most of them, due to Becker [2]), used for proving Theorems 1.1 and
1.2.
\par
\medskip
{\bf Lemma 2.4.} {\it Let $E$ be a hereditarily Pythagorean field
and $\sigma $ a generator of ${\cal G}(E(\sqrt{-1})/E)$. Then:}
\par
(i) $P(E) = \{2\}${\it , {\rm Br}$(E)$ has exponent $2$, and for
each integer $m \ge 3$, $E(\sqrt{-1})$ contains a primitive $m$-th
root of unity $\varepsilon _{m}$ and $\sigma (\varepsilon _{m}) =
\varepsilon _{m} ^{-1}$;}
\par
(ii) $E(\sqrt{-1}) ^{\ast } = E ^{\ast }E(\sqrt{-1}) ^{\ast n}${\it
, provided that $n \in \hbox{\Bbb N}$ and $2 \not\vert n$; moreover,
the natural embedding of $E$ in $E(\sqrt{-1})$ induces a group
isomorphism $E ^{\ast }/E ^{\ast n} \cong E(\sqrt{-1}) ^{\ast }/$
\par \noindent
$E(\sqrt{-1}) ^{\ast n}$;}
\par
(iii) $E ^{\ast }/E ^{\ast p ^{k}}$ {\it is presentable as a direct
product of $m$ isomorphic copies of $C _{p ^{k}}$, whenever $m, k
\in \hbox{\Bbb N}$ and $p \in \Pi _{m} (E)$.}
\par
\medskip
{\it Proof.} By Becker's theorem (cf. [2] and [3, (3.3)]), ${\cal G} _{E}$
is isomorphic to the topological semidirect product $G
_{E(\sqrt{-1})} \times \{\sigma \}$, where ${\cal G} _{E(\sqrt{-1})}$ is
abelian and $\sigma \tau \sigma = \tau ^{-1}$, for every $\tau \in G
_{E(\sqrt{-1})}$. Hence, by Galois theory, ${\cal G}(E _{\rm Ab}/E)$ is a
group of exponent $2$, where $E _{\rm Ab}$ is the maximal abelian
extension of $E$ in $\overline E$, and since $\sqrt{-1} \not\in E$,
the assertion that $P(E) = \{2\}$ becomes obvious. Observe now that
the field $\hbox{\Bbb Q} (\varepsilon _{m} + \varepsilon _{m}
^{-1})$ is formally real and normal over $\hbox{\Bbb Q}$, and
$[\hbox{\Bbb Q} (\varepsilon _{m})\colon \hbox{\Bbb Q} (\varepsilon
_{m} + \varepsilon _{m} ^{-1})] = 2$, for each $m \in \hbox{\Bbb
N}$. Note also that $\hbox{\Bbb Q} (\varepsilon _{m})$ is nonreal
whenever $m \ge 3$. In view of the Artin-Schreier theory and the
normality of $\hbox{\Bbb Q} (\varepsilon _{m})/\hbox{\Bbb Q}$, this
means that $\hbox{\Bbb Q} (\varepsilon _{m}) = \hbox{\Bbb Q}
(\varepsilon _{m} + \varepsilon _{m} ^{-1}) (\sqrt{-\beta _{m}})$,
for some $\beta _{m} \in \hbox{\Bbb Q} (\varepsilon _{m} +
\varepsilon _{m} ^{-1})$ presentable as a sum of elements of
$\hbox{\Bbb Q} (\varepsilon _{m} + \varepsilon _{m} ^{-1}) ^{\ast
2}$ (see [12, Ch. XI, Proposition 2 and the example in Sect. 2]). It
is therefore clear from the Pythagorean property of $E$ that $\beta
_{m} \in E$ and $E(\sqrt{-1}) = E(\varepsilon _{m})$, $m > 2$. As
${\cal G}(E _{\rm Ab}/E)$ has exponent $2$, this enables one to deduce from
the Merkurjev-Suslin theorem (cf. [14, (16.3) and (16.6)] and e.g.
[7, Lemma 3.3]) that Br$(E) _{2}$ has exponent $2$ and Br$(E) _{p} =
\{0\}$, for every $p \in \overline P \setminus \{2\}$. The obtained
results prove Lemma 2.4 (i). For the proof of the former assertion
of Lemma 2.4 (ii), it suffices to consider the special case where $n
\in \overline P \setminus \{2\}$. Since $n \not\in P(E)$, i.e. $E$
has no cyclic extensions of degree $n$, it follows from Lemma 2.4
(i) and [1, Ch. IX, Theorem 15] that $\sigma (\lambda )\lambda ^{-1}
\in E ^{\ast } \cap E(\sqrt{-1}) ^{\ast n}$, for all $\lambda \in
E(\sqrt{-1}) ^{\ast }$, which proves the former part of Lemma 2.4
(ii). The latter assertion of Lemma 2.4 (ii) is implied by the
former one, so it remains for us to prove Lemma 2.4 (iii). It
follows from Galois theory, the commutativity of ${\cal G} _{E(\sqrt{-1})}$
and Lemma 2.4 (i) that the character group of ${\cal G}(E(\sqrt{-1})
(p)/E(\sqrt{-1}))$ is a nontrivial divisible abelian torsion group,
for each $p \in \Pi (E)$. In view of [10, Theorem 23.1] and
Pontrjagin's duality theory, this amounts to saying that
${\cal G}(E(\sqrt{-1}) (p)/E(\sqrt{-1}))$ is presentable as a topological
product of isomorphic copies of $\hbox{\Bbb Z} _{p}$. It is
therefore clear from [4, (1.2)] and [17, Ch. I, Proposition 14] that
$p \in \Pi _{m} (E)$ if and only if ${\cal G}(E(\sqrt{-1})
(p)/E(\sqrt{-1})) \cong \hbox{\Bbb Z} _{p} ^{m}$. As $\varepsilon
_{p ^{k}} \in E(\sqrt{-1})$ when $k \in \hbox{\Bbb N}$, this enables
one to deduce Lemma 2.4 (iii) from Kummer theory and Lemma 2.4 (ii).
\par
\medskip
{\bf Remark 2.5.} In the setting of Lemma 2.4, with its proof, one
obtains from [3, (3.3)] that if $E$ has a unique ordering, then
${\rm cd}_{2} ({\cal G} _{E(\sqrt{-1})}) = 0$ and $E _{\rm Ab} = 
E(\sqrt{-1})$, whence Br$(L) _{2} = \{0\}$, for every $L \in {\rm
Fe}(E(\sqrt{-1}))$.
\par
\vskip0.6truecm \centerline{\bf 3. Proofs of the main results}
\par
\medskip
Let us note that, by [4, Proposition 3.1], a formally real field $E$
is quasilocal if and only if it is hereditarily Pythagorean with a
unique ordering, and cd$_{p} ({\cal G} _{E}) \le 2$, $p \in
P(E(\sqrt{-1}))$. This occurs if and only if $2 \not\in \Pi
(E(\sqrt{-1}))$ and ${\cal G} _{E}$ is isomorphic to the topological
semidirect product ${\cal G} _{E(\sqrt{-1})} \times \langle \sigma \rangle
$, where ${\cal G} _{E(\sqrt{-1})}$ is isomorphic to the topological group
product $\prod _{p \in P(E(\sqrt{-1}))} \hbox{\Bbb Z} _{p} ^{c(p)}$,
$c(p) = {\rm cd}_{p}({\cal G} _{E})$, $p \in P(E(\sqrt{-1}))$, $\sigma ^{2}
= 1$, and $\sigma \tau \sigma ^{-1} = \tau ^{-1}\colon \ \tau \in G
_{E(\sqrt{-1})}$ (cf. [2, Theorem 1], [3, (3.3)] or [4, (1.2) and
Proposition 3.1]). The following result proves the existence part of
Theorem 1.1 (i), and contains an analogue to the local reciprocity
law.
\par
\medskip
{\bf Proposition 3.1.} {\it Let $E$ be a hereditarily Pythagorean
field with a unique ordering, and let $R$ and $R _{1}$ be elements 
of {\rm Fe}$(E)$, such that $R \neq E$ and $R _{1} = R(\sqrt{-1})$. 
Assume that $R _{0}$ is a maximal subfield of $R$ with respect to 
the property that $2 \not\vert [R _{0}\colon E]$, and put $E _{1} = 
E(\sqrt{-1})$. Then:}
\par
(i) $R$ {\it equals $R _{0}$ or $R _{1}$, and the latter occurs if
and only if $R/E$ is normal;}
\par
(ii) $R/E$ {\it possesses an intermediate field $R ^{\prime }$
satisfying the following:}
\par
($\alpha $) $N(R ^{\prime }/E) = N(R/E)$ {\it and the prime divisors
of $[R ^{\prime }\colon E]$ and $[R\colon R ^{\prime }]$ lie in $\Pi
_{2} (E) \cup \{2\}$ and $\Pi _{1} (E)$, respectively;}
\par
($\beta $) {\it If $E$ is quasilocal, then $E ^{\ast }/N(R/E) \cong
{\cal G}(R ^{\prime } (\sqrt{-1})/E _{1}) \times C _{\delta }$, where
$\delta = [R\colon R _{0}]$; in particular, $R ^{\prime }$ is a
class field of $E$ associated with $N(R/E)$.}
\par
\medskip
{\it Proof.} Let $M \in {\rm Fe}(E)$ be a normal extension of $E$
including $R _{1}$. It follows from the structure of ${\cal G} _{E}$ that
$2 \vert [M\colon E]$, $4 \not\vert [M\colon E]$ and ${\cal G}(M/E)$ has a
subgroup $H$ of order $o(H) = [M\colon E]/2$. It is also clear that
$\varphi h\varphi ^{-1} = h ^{-1}\colon \ h \in H$, for every
$\varphi \in {\cal G}(M/E)$ of order $2$. Observing that $H$ is abelian and
normal in ${\cal G}(M/E)$, one obtains that subgroups of ${\cal G}(M/E)$ of odd
orders are included in $H$ and are normal in ${\cal G}(M/E)$, whereas
every $H _{2} \le {\cal G}(M/E)$ of even order equals its normalizer in
${\cal G}(M/E)$. Our argument also indicates that if $H _{0} \le {\cal G}(M/E)$
and $n \in \hbox{\Bbb N}$ divides $[M\colon E]$ and is divisible by
$o(H _{0})$, then there exists $H _{1} \le {\cal G}(M/E)$, such that $o(H
_{1}) = n$ and $H _{0} \subseteq H _{1}$. These results enable one
to deduce Proposition 3.1 (i) from Galois theory, and to prove that
$E$ has an extension $R ^{\prime }$ in $R$ with $[R ^{\prime }\colon
E]$ and $[R\colon R ^{\prime }]$ satisfying condition ($\alpha $) of
Proposition 3.1 (ii). Hence, by Lemma 2.3, $N(R/R ^{\prime }) = R
^{\prime \ast }$, which yields $N(R/E) = N(R ^{\prime }/E)$ as well.
It remains for us to prove Proposition 3.1 (ii) ($\beta $), so we
assume further that $E$ is quasilocal. Clearly, $R ^{\prime }$, $R
^{\prime } \cap R _{0} = R _{0} ^{\prime }$ and $R _{1} ^{\prime } =
R _{0} ^{\prime } (\sqrt{-1})$ are related in the same way as $R$,
$R _{0}$ and $R _{1}$. Taking also into account that Br$(E _{1})
_{p}$ is isomorphic to the quasicyclic $p$-group $\hbox{\Bbb Z} (p
^{\infty })$, for each $p \in \overline P$ dividing $[R _{0}
^{\prime }\colon E]$ (cf. [7, Theorem 3.1], [4, (1.2)], and [14,
(11.5)]), one obtains from [8, Theorem 3.1] that $E _{1}$ admits
local $p$-class field theory (i.e. subfields of $E (p)$ lying in
$\Omega (E)$ are related as in (2.1)), for each $p \in \overline P$
dividing $[R _{0} ^{\prime }\colon E]$. These results, combined with
[8, Theorem 3.1] and Lemma 2.2 indicate that $E _{1} ^{\ast }/N(R
_{1} ^{\prime }/E _{1}) \cong {\cal G}(R _{1} ^{\prime }/E _{1})$. Note
also that Lemmas 2.1, 2.3 and 2.4 (ii) lead to the following
conclusion:
\par
\medskip
(3.1) $N(R _{0}/E) = E ^{\ast } \cap N(R _{1}/E _{1})$, $N(R _{1}/E
_{1}) = N(R _{0}/E _{0})E _{1} ^{\ast \rho }$, where $\rho = [R
_{0}\colon E]$, and the embedding of $R _{0} ^{\prime }$ into $R
_{1} ^{\prime }$ induces a group isomorphism $E ^{\ast }/N(R _{0}/E)
\cong $ $E _{1} ^{\ast }/N(R _{1}/E _{1})$. Also, the set $I _{1}$ of
subgroups of $E _{1} ^{\ast }$ including $N(R _{1}/E _{1})$ consists
of the norm groups of the extensions of $E _{1}$ in $R _{1} ^{\prime
}$, and the mapping of $I _{1}$ into the set $I _{0}$ of subgroups
of $E ^{\ast }$ including $N(R _{0}/E)$, defined by the rule $H _{1}
\to H _{1} \cap E\colon \ H _{1} \in I _{1}$, is bijective. In
particular, $I _{0}$ coincides with the set of norm groups of
extensions of $E$ in $R _{0} ^{\prime }$.
\par
\medskip
\noindent
Observing finally that $E ^{\ast }/N(R _{1}/E) \cong E ^{\ast
}/N(R _{0}/E) \times E ^{\ast }/N(E _{1}/E)$, one completes the
proof of Proposition 3.1 (ii).
\par
\medskip
{\it Proofs of Theorems 1.1 and 1.2.} Assume that $E$, $E _{1}$,
$\Pi _{1} (E)$ and $\Pi _{2} (E)$ satisfy the conditions of
Proposition 3.1, fix a field $\Psi \in {\rm Fe}(E)$ and a group $U
\in {\rm Nr}(E)$, take class fields $R$ and $T$ of $E$ associated
with $U$, and as in Proposition 3.1, attach fields $R _{0}$, $R
_{1}$, $T _{0}$, $T _{1}$ and $\Psi _{0}$, $\Psi _{1}$ to $R$, $T$
and $\Psi $, respectively. Applying Lemma 2.1 and the former two
assertions of (3.1), one obtains consecutively that $N(R _{0}/E) =
N(T _{0}/E)$, $N(R _{1}/E _{1}) = N(T _{1}/E _{1})$ and $N(R _{1}/E)
= N(T _{1}/E)$. Hence, by Proposition 3.1 (ii) ($\alpha $) and [8,
Theorem 3.1], we have $R _{1} = T _{1}$. Using Lemma 2.1 and [8,
Theorem 3.1], one also sees that $N(\Psi /E) \subseteq U$ if and
only if $N(\Psi _{1}/E _{1}) \subseteq N(R _{1}/E _{1})$, which
holds if and only if $R _{1} \subseteq \Psi _{1}$. Furthermore, $R
_{0}$ and $T _{0}$ are the fixed fields of some Sylow $2$-subgroups
of ${\cal G}(R _{1}/E)$, and therefore, are $E$-isomorphic. As $R _{0}$ is
a maximal subfield of $R _{1}$, it becomes clear that $R _{0}T _{0}
= R _{1}$, unless $R _{0} = T _{0}$. It is now easily deduced from
Galois theory and Sylow's theorems that $R _{1} \subseteq \Psi _{1}$
if and only if $R _{0}$ embeds in $\Psi _{0}$ over $E$. These
observations, combined with Proposition 3.1, prove Theorem 1.1 (i)
and (ii).
\par
We turn to the proof of Theorem 1.2 (ii). When $2 \not\vert n$, the
equality $E ^{\ast n} = E ^{\ast n(E)}$ is obvious. Since, by [7,
Lemma 3.5], $E (2) = E(\sqrt{-1})$, and by Kummer theory, this
yields $E ^{\ast } = E ^{\ast 2} \cup (-1)E ^{\ast 2}$, one also
sees that $E ^{\ast 2 ^{k}n} = E ^{\ast 2n(E)}$ and $E ^{\ast }/E
^{\ast 2 ^{k}n} \cong [E ^{\ast }/E ^{\ast n(E)}] \times C _{2}$,
for every $k \in \hbox{\Bbb N}$. At the same time, it follows from
Lemma 2.4 (iii) and the inequalities cd$_{p} ({\cal G} _{E}) \le 2$, $p \in
\overline P \setminus \{2\}$, that $E ^{\ast }/E ^{\ast n(E)} \cong 
$
\par \noindent
$C _{n(E)} \times C _{n _{1} (E)}$. Summing up the obtained results,
one proves Theorem 1.2 (ii) in general.
\par
Our next objective is to prove Theorems 1.1 (iii). Assume that $n =
n _{1} (E)$, i.e. $E ^{\ast }/E ^{\ast n} \cong C _{n} \times C
_{n}$, fix elements $b _{1}$ and $b _{2}$ of $E ^{\ast }$ so that
$\langle b _{j}E ^{\ast n}: j = 1, 2\rangle = $ $E ^{\ast }/E ^{\ast
n}$, and denote by $\Phi _{1}$ the extension of $E _{1}$ in
$\overline E$ obtained by adjoining the $n$-th roots of $b _{1}$ and
$b _{2}$. Then $\Phi _{1}/E _{1}$ is clearly a Kummer extension with 
${\cal G}(\Phi _{1}/E _{1}) \cong C _{n} \times C _{n} \cong E _{1} ^{\ast 
}/E _{1} ^{\ast n}$. Since $\Pi (E) \setminus \{2\} = \Pi _{1} (E) 
\cup \Pi _{2} (E)$, this implies in conjunction with [8, Theorem 
3.1], Proposition 3.1 and Lemmas 2.2 and 2.4 that $N(\Phi _{1}/E 
_{1}) = E _{1} ^{\ast n}$, $N(\Phi _{1}/E) = E ^{\ast 2n}$ and 
$N(\Phi _{0}/E) = E ^{\ast n}$, where $\Phi _{0}$ is a real closure 
of $E$ in $\Phi _{1}$. As $E ^{\ast n}$ is included in every 
subgroup of $E ^{\ast }$ of index $n$, the obtained result enables 
one to deduce Theorem 1.1 (iii) from Proposition 3.1 and statement 
(3.1).
\par
We are now in a position to complete the proof of Theorem 1.2. The
normality of $\Lambda _{2}/E$ follows at once from Proposition 3.1
and Galois theory, the conclusion of Theorem 1.2 (iii) is contained
in Proposition 3.1, and the uniqueness of $\Lambda _{0}$ (up-to an
$E$-isomorphism) follows from the prosolvability of ${\cal G} _{E}$ and the
fact that $\Lambda _{0}$ is the fixed field of a Hall pro-$(\Pi _{1}
(E) \cup \{2\})$-subgroup of ${\cal G} _{E}$. What remains to be seen is
the validity of Theorem 1.2 (i). Let $L \subseteq \overline E$ be a
class field of $E$, such that $2 \not\vert [L\colon E]$. In view of
the prosolvability of ${\cal G} _{E}$, Theorem 1.1 (iii) and Galois theory,
then ${\cal G} _{L}$ includes some Hall pro-$(\Pi _{1} (E) \cup
\{2\})$-subgroup of ${\cal G} _{E}$, which ensures that $L$ embeds in
$\Lambda _{0}$ over $E$. Returning to the proof of Theorem 1.1 (i)
and (ii), one immediately obtains that the $E$-isomorphic copy of
$L$ in $\Lambda _{0}$ is unique. It is therefore easy to deduce that
$\Sigma _{0} \cup \Sigma _{1}$ is a system of representatives of the
class fields of $E$, and that the canonical mapping of $\Sigma _{0}
\cup \Sigma _{1}$ into Nr$(E)$ has the properties required by
Theorem 1.2 (i). Thus Theorems 1.1 and 1.2 are proved.
\par
\medskip
{\bf Remark 3.2.} Fix $E$, $\Pi _{1} (E)$ and $\Pi _{2} (E)$ as in
Proposition 3.1, and suppose that $E$ is quasilocal.
\par
(i) Arguing as in the proof of Theorem 1.2 (ii), one obtains that if
$H$ is a subgroup of $E ^{\ast }$ of index $n \in \hbox{\Bbb N}$,
and $E ^{\ast }/H$ is of exponent $e$, then $n \vert e ^{2}$, $n/e$
is not divisible by any $p \in \Pi _{1} (E) \cup \{2\}$, and $E
^{\ast }/H \cong C _{e} \times C _{n/e}$.
\par
(ii) The proof of Proposition 3.1 shows that if $R$ is a nonreal
finite extension of $E$, then Br$(R)$ is isomorphic to the direct
sum $\oplus _{p \in \Pi _{2} (E)} \hbox{\Bbb Z} (p ^{\infty })$.
Hence, by Lemma 2.4 (i), $E$ is SQL if and only if $\Pi _{1} (E) =
\phi $.
\par
(iii) It follows from Proposition 3.1 and Remark 2.5 that if $R \in
{\rm Fe}(E)$, then $R \cap E _{\rm Ab} \subseteq {\rm cl}(R/E)$. It
is also clear that cl$(R/E) = R \cap E _{\rm Ab}$ if and only if
$[R\colon E]$ is not divisible by any $p \in \Pi _{2} (E)$. This
holds if and only if $R$ embeds over $E$ into the fixed field of
some Hall pro-$\Pi _{2} (E)$-subgroup of ${\cal G} _{E}$.
\par
\medskip
Note finally that if $P _{1} \subset P \subseteq \overline P$ and $2
\in P \setminus P _{1}$, then there exists a formally real
quasilocal field $F$ with $\Pi (F) = P$ and $\Pi _{1} (F) = P _{1}$.
Also, it follows from [4, Proposition 3.1] that ${\cal G} _{F}$ does not
depend on the choice of $F$ but is uniquely determined by $(P, P
_{1})$, up-to an isomorphism.
\vskip1cm \centerline{ REFERENCES} \vglue15pt\baselineskip12.8pt
\def\num#1{\smallskip\item{\hbox to\parindent{\enskip [#1]\hfill}}}
\parindent=1.38cm
\par
\medskip
\num{1} A.A. {\pc ALBERT}, {\sl Modern Higher Algebra.} Univ. of Chicago
Press, XIV, Chicago, Ill., 1937.
\par
\num{2} E. {\pc BECKER}, {\sl Hereditarily-pythagorean Fields and Orderings of
higher Level.} IMPA Lecture Notes {\bf 29}, IMPA, Rio de Janeiro,
1978.
\par
\num{3} L. {\pc BR\"OCKER}, {\sl Characterization of fans and hereditarily
Pythagorean fields.} Math. Z. {\bf 151} (1976), 149-163.
\par
\num{4} I.D. {\pc CHIPCHAKOV}, {\sl On the Galois cohomological dimensions of
stable fields with Henselian valuations.} Comm. Algebra {\bf 30}
(2002), 1549-1574.
\par
\num{5} I.D. {\pc CHIPCHAKOV}, {\sl Class field theory for strictly
quasilocal fields with Henselian discrete valuations.} Manuscr. math.
{\bf 119} (2006), 383-394.
\par
\num{6} I.D. {\pc CHIPCHAKOV}, {\sl On nilpotent Galois groups and the scope
of the norm limitation theorem in one-dimensional abstract local
class field theory.} In: Proc. of ICTAMI 05, Alba Iulia, Romania,
15.9-18.9. 2005; Acta Univ. Apulensis {\bf 10} (2005), 149-167.
\par
\num{7} I.D. {\pc CHIPCHAKOV}, {\sl On the residue fields of Henselian valued
stable fields.} J. Algebra {\bf 319} (2008), No1, 16-49.
\par
\num{8} I.D. {\pc CHIPCHAKOV}, {\sl One-dimensional abstract local class
field theory.} E-print, arXiv:math/0506515v4 [math.RA].
\par
\num{9} I.D. {\pc CHIPCHAKOV}, {\sl On the Brauer groups of quasilocal
nonreal fields and the norm groups of their finite Galois
extensions.} E-print, arXiv:math/
\par \noindent
0707.4245v3 [math.RA].
\par
\num{10} L. {\pc FUCHS}, {\sl Infinite Abelian Groups.} V. I and II, Pure and
Applied Mathematics, Academic Press, IX and XI, New York-London,
1970.
\par
\num{11} K. {\pc IWASAWA}, {\sl  Local Class Field Theory.} Iwanami Shoten,
Japan, 1980 (Japanese: English transl. in Oxford Mathematical
Monographs. Oxford Univ. Press, New York; Clarendon Press, VIII,
Oxford, 1986).
\par
\num{12} S. {\pc LANG}, {\sl Algebra.} Addison-Wesley Publ. Comp., Reading,
Mass., 1965.
\par
\num{13} A.S. {\pc MERKURJEV}, {\sl Brauer groups of  fields.}  Comm. Algebra
{\bf 11} (1983), No 22, 2611-2624.
\par
\num{14} A.S. {\pc MERKURJEV}, A.A. {\pc SUSLIN}, {\sl $K$-cohomology  of
Severi-Brauer varieties and the norm residue homomorphism.} Izv.
Akad. Nauk SSSR {\bf 46}, No 5 (1982), 1011-1046 (in Russian:
English transl. in Math. USSR Izv. {\bf 21} (1983), 307-340).
\par
\num{15} R. {\pc PIERCE}, {\sl Associative Algebras.} Graduate Texts in Math.
{\bf 88}, Springer-Verlag, New York-Heidelberg-Berlin, 1982.
\par
\num{16} L.S. {\pc PONTRJAGIN}, {\sl Topological Groups, 4th Ed.} Nauka,
Moscow, 1984 (in Russian).
\par
\num{17} J.-P. {\pc SERRE}, {\sl Cohomologie Galoisienne.} Cours au College de
France, 1962-1963. Lecture Notes in Math. {\bf 5}, Springer-Verlag,
Berlin-Heidelberg-New York, 1965.
\par
\medskip
Institute of Mathematics and Informatics, Bulgarian Academy \par
of Sciences, Acad. G. Bonchev Str., bl. 8, Sofia, 1113, Bulgaria
\par
E-mail chipchak@math.bas.bg
\vskip1cm
\def\pc#1{\eightrm#1\sixrm}
\hfill\vtop{\eightrm\hbox to 5cm{\hfill Ivan {\pc
CHIPCHAKOV}\hfill}
 \hbox to 5cm{\hfill Institute of Mathematics and Informatics\hfill}\vskip-2pt
 \hbox to 5cm{\hfill Bulgarian Academy of Sciences\hfill}
\hbox to 5cm{\hfill Acad. G. Bonchev Str., bl. 8\hfill} \hbox to
5cm{\hfill 1113 {\pc SOFIA,} Bulgaria\hfill}}
\end